\documentclass{amsart}

\newtheorem{theorem}{Theorem}[section]
\newtheorem{lemma}[theorem]{Lemma}
\newtheorem{prop}[theorem]{Proposition}

\theoremstyle{definition}

\def\la{{\langle}}
\def\ra{{\rangle}}
\def\pr{{\bf pr}} 

\def\a{{\mathfrak a}}
\def\b{{\mathfrak b}}
\def\g{{\mathfrak g}}
\def\k{{\mathfrak k}}
\def\n{{\mathfrak n}}
\def\m{{\mathfrak m}}
\def\p{{\mathfrak p}}
\def\t{{\mathfrak t}}
\def\z{{\mathfrak z}}

\def\C{{\mathbb C}}
\def\R{{\mathbb R}}

\theoremstyle{remark}

\numberwithin{equation}{section}
\newcommand{\abs}[1]{\lvert#1\rvert} 
\begin{document}
\author[B.Kr\"otz]{Bernhard Kr\"otz}
\author[M. Otto]{Michael Otto}
\address{Department of Mathematics, 1222 
University of Oregon, Eugene OR 97403-1222}
\address{The Ohio State University, Department of Mathematics, 
231 West 18th Avenue, 
Columbus, OH 43210-1174}
\email{kroetz@math.uoregon.edu}
\email{otto@math.ohio-state.edu}
\title[Complex convexity theorem]
{A refinement of the complex convexity theorem via 
symplectic techniques}
\thanks{The first author was supported in part by NSF grant DMS-0097314}

\begin{abstract} 
We apply techniques from symplectic geometry to extend and give a new proof
of the complex convexity theorem of Gindikin-Kr\"otz.
\end {abstract}
\maketitle
\section{Introduction}

Let $G$ be a semisimple Lie group and $G=NAK$ 
be an Iwasawa decomposition. Let $X\in \a={\rm Lie}(A)$ be such that 
${\rm Spec}({\rm ad} X)\subseteq ]-{\pi\over 2}, {\pi\over 2}[$ and set
$a=\exp(iX)$. It is known that $Ka\subseteq N_\C A_\C K_\C$. 
The main result of this paper asserts that 
\begin{equation}\label{1000}
\Im \log \tilde  a(Ka)={\rm conv}({\mathcal W}.X)\ .
\end{equation}
Here $\tilde a: N_\C A_\C K_\C\to A_\C$ denotes a middle projection 
and ${\rm conv}({\mathcal W}.X)$ stands for the convex hull 
of the Weyl group orbit ${\mathcal W}.X$.  
We note that the inclusion $''\subseteq''$ in 
(\ref{1000}) is the complex convexity theorem from 
\cite{GK}. 

\par We prove (\ref{1000}) with symplectic methods. As a byproduct
we obtain interesting new classes of compact Hamiltonian  manifolds 
and Lagrangian submanifolds thereof.

\par Let us now be more specific about the used techniques. 
We first consider the case where $G$ is a complex group. In  this 
situation we show that $B_\C=A_\C N_\C$ carries a natural 
structure of a Poisson Lie group. Locally, we can identify 
$B_\C$ inside of $G_\C/K_\C$ and consequently we obtain 
a local action of $G_\C$ on $B_\C$. Within this identification 
the symplectic leaf 
$P_a$ through $a=\exp(iX)\in B_\C$ becomes a local $K_\C$-orbit. 
Interestingly, the symplectic form on $P_a$ remains non-degenerate on the 
totally real $K$-orbit  $M_a\subseteq P_a$.
We then exhibit a Hamiltonian torus action on the compact 
symplectic manifold $M_a$ and show that that (\ref{1000}) becomes a consequence of the 
Atiyah-Guillemin-Sternberg convexity theorem. 
Finally, 
the general case of arbitrary $G$ can be handled by descent 
to certain Lagrangian submanifolds $Q_a\subseteq M_a$ by means 
of the recently discovered convexity theorem
\cite{KO}.

\section{Notation and basic facts} 

In this section we recall some basic facts about semisimple Lie algebras and
groups. We make an emphasis on the complexified Iwasawa decomposition. 
Furthermore we review some standard facts on the double complexification 
of a semisimple Lie algebra. 

\subsection{Iwasawa decomposition and complex crown}
We let $\g$ denote a semisimple Lie algebra and let 
$\g=\k+\p$ be  a Cartan decomposition of $\g$.
For a maximal abelian subspace $\mathfrak{a}$
of $\mathfrak{p}$ let $\ \Sigma = \Sigma (\mathfrak{g},
\mathfrak{a} ) \subseteq \mathfrak{a}^{*} \ $ be the 
corresponding root system. 
Then $\mathfrak{g}$ admits a root space decomposition
$$ \mathfrak{g} = \mathfrak{a} \oplus \mathfrak{m} \oplus
\bigoplus_{\alpha \in \Sigma} \mathfrak{g}^{\alpha} ,$$
where $\ \mathfrak{m}=\mathfrak{z}_{\mathfrak{k}} (\mathfrak{a})
\ $ and $\ \mathfrak{g}^{\alpha}=\{ X \in \mathfrak{g} :
(\forall H\in \mathfrak{a}) \  \left[ H,X \right] = \alpha(H)X\} $. \\
For a fixed positive system $\Sigma^{+}$ define $\
\mathfrak{n} = \bigoplus_{\alpha \in \Sigma^{+}}
\mathfrak{g}^{\alpha} \ $ . Then we have the Iwasawa 
decomposition on the Lie algebra level:
\[ \mathfrak{g} = \mathfrak{n} \oplus \mathfrak{a} \oplus
\mathfrak{k} . \]

\par For any real Lie algebra $\mathfrak{l}$ we write 
$\mathfrak{l}_{\mathbb C}$ for its complexification. \\
\par In the following $G_{\mathbb C}$ will denote  a simply connected Lie group with Lie
algebra $\mathfrak{g}_{\mathbb C}$. We write  $\ G, K, K_{\mathbb C}, 
A, A_{\mathbb C}, N \ $ and $N_{\mathbb C}$ for the analytic 
subgroups of $G_{\mathbb C}$ corresponding to the subalgebras $\ \mathfrak{g}, 
\mathfrak{k}, \mathfrak{k}_{\mathbb C}, \mathfrak{a}, 
\mathfrak{a}_{\mathbb C}, \mathfrak{n} \ $ and 
$\mathfrak{n}_{\mathbb C}$ , respectively. \\
The Weyl group of $\Sigma$ can be defined 
by $\mathcal{ W}=N_{K}(\mathfrak{a})/Z_{K}(\mathfrak{a})
$.

Following \cite{AG} we define a bounded and convex subset of $\mathfrak{a}$ that plays a central role:
\begin{equation} \Omega = \{ X\in \mathfrak{a} : \abs{\alpha(X)}<
\frac{\pi}{2} \ \forall \alpha \in \Sigma \} . 
\end{equation}
\\
With $\Omega$ we define a $G-K_{\mathbb C}$-double coset domain in 
$G_{\mathbb C}$ by 
$$ \tilde \Xi=G\exp(i\Omega)K_{\mathbb C}\ .$$
Also we write 
$$\Xi=\tilde \Xi/ K_{\mathbb C}$$
for the union of right $K_{\mathbb C}$-cosets of $\tilde \Xi$ in the complex symmetric space 
$G_{\mathbb C}/ K_{\mathbb C}$. We refer to $\Xi$ as the {\it complex crown}
of the symmetric space $G/K$. Notice that $\Xi$ is independent of the 
choice of $\mathfrak a$, hence generically defined through $G/K$. \\
It is known that $\tilde{\Xi}$ is an open and $G$-invariant
subset of $N_\C A_\C K_\C$ (cf.\ \cite{M} for a short 
proof).  

\par Next we consider the open and dense cell 
$N_\C A_\C K_\C\subseteq G_\C$ in more detail. Define 
$F=K_\C \cap A_\C$ and recall that $F=K\cap\exp(i\a)$ is a 
finite $2$-group. 
Standard techniques imply that the mapping 
\begin{equation}\label{200} 
N_\C \times [A_\C\times_F K_\C]\to N_\C A_\C K_\C, \ \ (n,[a,k])
\mapsto nak
\end{equation}
is a biholomorphism. In particular (\ref{200}) induces 
a holomorphic map $\tilde n:  N_\C A_\C K_\C\to N_\C$ 
and a multi-valued holomorphic mapping $\tilde a:  N_\C A_\C K_\C\to A_\C$ such that 
$x\in \tilde n(x)\tilde a(x) K_\C$ for all $x\in  N_\C A_\C K_\C$. 

\par Set $B_\C = N_\C A_\C$. Clearly, 
$\b_\C=\a_\C+\n_\C$ is the Lie algebra of $B_\C$.
We define a multi-valued holomorphic  map $\tilde b: N_\C A_\C K_\C \to B_\C$ by 
$\tilde b(x)=\tilde n(x) \tilde a(x)$. 

\par Recall that $\Xi$ is contractible, and hence 
simply connected. It follows that 
the restriction $\tilde a|_{\tilde \Xi}$ has 
a unique single-valued holomorphic lift $\log \tilde a: \tilde \Xi\to \a_\C$
such that $\log \tilde a({\bf 1})=0$. 
Consequently, $\tilde b|_{\tilde \Xi}$ lifts to a single-valued holomorphic 
map $\tilde\Xi \to B _\C$ which shall also be denoted by $\tilde b$. 

\subsection{Double complexification of a Lie algebra} For the remainder of this 
section we will assume that $\g$ carries a complex structure, say $j$. Then $\g$  can be viewed as the 
complexification of its compact real form $\k$ and the Cartan decomposition becomes 
$\g=\k+j\k$. The Cartan involution  $\theta$ on $\g$ coincides 
with the complex conjugation $\bar{}:\g\rightarrow \g$ with respect 
to $\k$. \\
A second complexification yields $\g_\C$, which carries another complex
structure $i:\g_\C \rightarrow \g_\C$. \\
The following map $\varphi$ defines a real Lie algebra isomorphism.
\[ \varphi :\g_\C \rightarrow \g\times \g, \quad X+iY \mapsto 
(X+jY,\bar{X}+j\bar{Y}) . \]
Under $\varphi$ we have the following identifications:
\begin{eqnarray}
\k_\C & = & \{ (Z,Z):Z\in \g \} , \\
\a_\C & = & \{ (Z,-Z):Z\in \a+j\a \} , \\
\n_\C & = & \{ (X+jY,\theta(X-jY)):X,Y\in \n \} = 
\n \times \bar{\n} ,
\end{eqnarray}
where $\bar{\n}=\bigoplus_{-\alpha \in \Sigma^{+}} \g^\alpha $.

\section{The complex convexity theorem for complex groups}
The objective of this section is to provide a symplectic proof
of the complex convexity theorem in the case of $G$ complex. 
If $G$ is complex, then we can endow $B_\C$ with a natural 
structure of a Poisson Lie group. The symplectic leaves 
become local $K_\C$-orbits. We show that the totally real 
$K$-orbit in each leaf is again a symplectic manifold.  With that 
we obtain compact symplectic manifolds with appropriate 
Hamiltonian torus actions. The complex convexity theorem then becomes 
a consequence of the Atiyah-Guillemin-Sternberg convexity theorem. 

\subsection{The Poisson Lie group $B_\C$} Throughout 
this section $\g$ denotes a complex semisimple 
Lie algebra. 
Our first task is to define a bilinear form $\la , \rangle $ 
on $\g_\C = \g\times \g$ which gives $(\g_\C, \b_\C, \k_\C)$ 
the structure of a Manin triple. 
To that end let $\kappa$ be the Killing form of the 
complex Lie algebra $\g$. Then define on $\g_\C=\g\times \g$ 
a bilinear form by 
\[ \langle(X,Y),(X',Y')\rangle = 
\Re \kappa(X,X')-\Re \kappa(Y,Y') \]
for $(X,Y), (X',Y')\in \g_\C$. 

\begin{prop}
The bilinear form $\la, \ra$ on $\g_\C=\g\times \g$ has the following properties: 
\item {{\rm (i)}} $\la, \ra$ is symmetric, non-degenerate and $G_\C$-invariant. 
\item {{\rm (ii)}} $\la \b_\C, \b_\C \ra =\{0\}$ and $\la \k_\C, \k_\C\ra =\{0\}$. 
\end{prop}
\begin{proof}
(i) is immediate from the definition. Moving on (ii) we observe  
that the relations $\langle\k_\C,\k_\C\rangle=\{0\}$ and $\langle\a_\C,\a_\C\rangle
=\{0\}$ are straightforward 
from
the identifications (2.3)-(2.5). 
Finally, the fact that root spaces $\g^\alpha$ and $\g^\beta$ 
are $\kappa$-orthogonal if $\alpha+\beta \neq 0 $ 
implies $\langle\n_\C,\a_\C+\n_\C\rangle=\{0\}$.
\end{proof} 
Notice that Proposition 3.1 just says that $\la, \ra$ turns  
$(\g_\C, \b_\C, \k_\C)$ into a Manin triple. Accordingly 
$B_\C$ becomes a Poisson Lie group whose symplectic leaves are 
local $K_\C$-orbits (cf. \cite{CP}). Here we shall only be interested 
in the leaves through points $a=\exp(iX)$ for $X\in \Omega$. Denote the 
leaf containing $a$ by $P_a$. 
Then 
$$P_a=\{ \tilde b(ka)\in B_\C  : \ k\in K_\C, ka\in N_\C A_\C K_\C \}_0\, ,  $$ 
where $\{\cdot\}_0$ refers to the connected component of $\{\cdot\}$
containing $a$. 
\par For an element $Z\in \k_\C$ we write $\tilde Z$ for the corresponding 
vector field on $P_a$, i.e. if $b=\tilde b(ka)\in P_a$, then 
$$\tilde Z_b={d\over dt}\Big|_{t=0} \tilde b(\exp(tZ)ka)\ .$$

\par Write $\pr_{\k_\C}: \g_\C \to \k_\C$ and $\pr_{\b_\C}: \g_\C \to \b_\C$
for the projections along $\b_\C$, resp. $\k_\C$.  

\par We notice that $T_b P_a=\{\tilde Z_b : Z\in \k_\C\}$. 
The symplectic form $\tilde \omega$ on $P_a$ is then given by 
\begin{equation}\label{syform}
\tilde{\omega}_b (\tilde{Y}_b, \tilde{Z}_b) = \langle \pr_{\k_\C}({\rm Ad}(b^{-1})Y),
{\rm Ad}(b^{-1})Z \rangle = -\langle \pr_{\b_\C}({\rm Ad}(b^{-1})Y),
{\rm Ad}(b^{-1})Z \rangle \ ,
\end{equation}
for $Y,Z\in\k_\C$ (cf. (\cite{CP}).  The second equality in (\ref{syform})
follows from Proposition 3.1.

\subsection{The totally real $K$-orbit in the symplectic leaf}  

\par Our interest is not so much with $P_a$ as it is 
with its totally real submanifold 
$$M_a=\tilde b (Ka)\ .$$
Then $\tilde \omega$ induces a closed 2-form $\omega$ on 
$M_a$ by $\omega=\tilde\omega|_{T M_a \times T M_a}$.
A priori it is not clear that $\omega$ is non-degenerate, i.e. that  
$(M_a, \omega)$ is a symplectic manifold. This will be shown now. 
We start with a simple algebraic fact. 

\begin{lemma} With respect to $\la,\ra$ one has $\k^\perp=\k_\C +\b$. 
\end{lemma}
\begin{proof} Because of the non-degeneracy of $\la , \ra $,  
it is sufficient to verify  that $\k_\C +\b\subseteq \k^\perp$. Clearly, $\k_\C \subseteq \k^{\perp}$
since $\k_\C$ is isotropic. As $\b=\a +\n$, it thus remains to show that 
$\la \k, \a \ra=\{0\}$ and $\la \k, \n\ra=\{0\}$. 
Now $\la \k, \a \ra=\{0\}$ follows from (2.3-4) and $\Re \kappa (\k, \a)=\{0\}$. 
Finally we show that $\la \k, \n \ra =\{0\}$.  
For that fix an arbitrary element $W=(\sum_{\alpha\in \Sigma^+} Y_\alpha, \sum_{\alpha\in \Sigma^+} \bar Y_\alpha)$ 
of $\n\subseteq \g\times\g$; here $Y_\alpha\in \g^\alpha$. 
Likewise let $U=(V+\sum_{\alpha\in \Sigma^+}(Z_\alpha+\bar Z_\alpha), V+
\sum_{\alpha\in \Sigma^+}(Z_\alpha+\bar Z_\alpha)) $ be an element of $\k\subseteq \g\times \g$; here 
$V\in \m$ and $Z_\alpha\in \g^\alpha$. 
\par As $\m\perp_\kappa \g^\alpha$ and  
$\g^\alpha\perp_\kappa\g^\beta$ for $\alpha+\beta\neq 0$, we obtain 

\begin{align*}
\langle U, W \rangle & = \sum_{\alpha\in \Sigma^+} \Re \kappa (Z_\alpha+\bar Z_\alpha, 
Y_\alpha) -\Re\kappa(Z_\alpha +\bar Z_\alpha , \bar Y_\alpha) \\ 
&=  \sum_{\alpha\in \Sigma^+} \Re \kappa (\bar Z_\alpha, 
Y_\alpha) -\Re\kappa(Z_\alpha , \bar Y_\alpha)=0\ .
\end{align*}
This concludes the proof of the lemma. 
\end{proof}
 
\begin{lemma}
At any point $b\in M_a$, the bilinear form 
$\omega_b: T_b M_a\times T_b M_a \to  \R$  is non-degenerate.
\end{lemma}
\begin{proof}
Let $b=\tilde{b}(ka)\in M_a$. Then $b=kak'$ for some $k\in K, k' \in K_\C $. 
Notice that $T_bM_a=\{ \tilde Y_b: Y\in \k\}$. 
\par Assume that there is an $U\in \k$ such that $\omega_b(\tilde U_b, \tilde Y_b)=0$
for all $Y\in \k$, i.e. 
\[ \langle \pr_{\k_\C}({\rm Ad}(b^{-1})U),{\rm Ad}(b^{-1})Y \rangle =0 \quad \forall \ Y\in \k.
\]
We have to show that $\tilde U_b=0$. Set $Z=\pr_{\k_\C}({\rm Ad}(b^{-1})U). $ Then $\rm {Ad}(b)Z\in \k^\perp$. 
Thus Lemma 3.2 implies that 
\begin{equation}\label{eins}
Z \in {\rm Ad}(b^{-1})(\k_\C + \b) . 
\end{equation}
On the other hand, by definition,
\begin{equation}\label{zwei}
Z \in {\rm Ad}(b^{-1})\k + \b_\C = {\rm Ad}(b^{-1})(\k+\b_\C) .
\end{equation}
From (\ref{eins}) and (\ref{zwei}) it follows that
\[ Z \in {\rm Ad}(b^{-1})(\k_\C + \b) \cap {\rm Ad}(b^{-1})(\k+\b_\C)
= {\rm Ad}(b^{-1})\g . \]
Moreover, since $Z$ is an image point of $\pr_{\k_\C}$,
\[ Z \in {\rm Ad}(b^{-1})\g \cap \k_\C , \]
i.e.
\[ {\rm Ad}(b)Z \in \g \cap {\rm Ad}(b)\k_\C . \]
As $b=kak'$ we now get
\[ {\rm Ad}(b)Z \in \g \cap {\rm Ad}(k) {\rm Ad}(a)\k_\C
={\rm Ad}(k)(\g \cap {\rm Ad}(a)\k_\C )\ .\]
Using standard techniques (see \cite{AG} or \cite{M}, Lemma 2), it follows from (2.1) that 
$$\g \cap {\rm Ad}(a)\k_\C = \z_{\k}(X)\ .$$
Hence ${\rm Ad}(b)Z\in{\rm Ad}(k)\z_\k(X)$. 
\par 
From
\[ {\rm Ad}(b^{-1}) U \in Z+\b_\C \]
we conclude
\[ U \in {\rm Ad}(b)Z +\b_\C \subseteq {\rm Ad}(k) \z_{\k}(X)+\b_\C . \]
Since $U$ was assumed to lie in $\k$, we finally get 
\[ U \in {\rm Ad} (k) \z_{\k}(X) . \]
But this just means that $\tilde{U}_b =0 $, concluding the proof that $\omega_b$
is non-degenerate.
\end{proof}

\subsection{Hamiltonian torus action on the totally real leaf $M_a$}

It follows from Lemma 3.3 that $(M_a, \omega)$ is a
compact symplectic manifold. Clearly, the torus $T=\exp(j\a)$ acts 
on $M_a$ as $K$ does. We wish to show that the action of $T$ on $M_a$ 
is Hamiltonian and identify the corresponding moment 
map $\Phi: M_a\to \t^*$, where $\t=j\a$ is the Lie algebra 
of $T$.  We will identify $\t^*$ with $\a$ via 
the linear isomorphism 
\begin{equation}\label{ident}
\a\to \t^*, \ \ Y\mapsto (Z\mapsto \la iY,Z\ra)\, .
\end{equation}
\par We shall write $\pr_{\a_\C}: \g_\C\to \a_\C$ for the projection along
$\k_\C+\n_\C$.

\begin{prop}
The action of the torus $T=\exp( j\a)$ on $M_a$ is Hamiltonian with
momentum map
\[ \Phi:M_a \rightarrow \t^* \simeq \a, \ \tilde{b}(ka) \mapsto \Im \log \tilde{a}
(ka) . \]
\end{prop}
\begin{proof}   We first show that 
$T=\exp(j\a)$ acts on $M_a$ by symplectomorphisms. For that 
we first notice that $T$ normalizes $B_\C$. Hence the 
action of $T$ on $M_a\subseteq B_\C$  is  given by conjugation, i.e.
$$T\times M_a\to M_a, \ \ (t,b)\mapsto t.b=tbt^{-1}\, .$$
Moreover, for each $t\in T$ the map ${\rm Ad}(t)$ commutes both with 
$\pr_{\k_\C}$ and $\pr_{\b_\C}$. Combining these facts,  it is
then straightforward from (3.1)
that $T$ acts indeed symplectically on $M_a$. 
\par Next we show that 
\begin{equation} \label{100}
\omega(\tilde Y, \tilde Z)=0 \qquad \forall Y,Z\in\t \ .
\end{equation}
Fix $b\in M_a$.  From the definition 
(\ref{syform}) we obtain that 
$$\omega_b(\tilde Y_b, \tilde Z_b)=\la \pr_{\b_\C} ({\rm Ad}(b^{-1}) Y), 
{\rm Ad}(b^{-1} Z)\ra\ .$$
Now for $Y,Z\in\t$ we have ${\rm Ad}(b^{-1})Y\in Y+\n_\C$ and 
${\rm Ad}(b^{-1})Z\in Z+\n_\C$. 
From Proposition 3.1 we know $\la \n_\C, \n_\C\ra=\{0\}$ and 
$\la \k_\C, \k_\C\ra =\{0\}$. Hence, to prove (\ref{100}) it suffices to show
$\la \t,\n_\C\ra=\{0\}$. \\
Let $ U=(jZ,jZ) \in \t $ and $ V=(\sum_{\alpha\in \Sigma^+} X_\alpha +j Y_\alpha ,
\sum_{\alpha\in \Sigma^+} \bar{X}_\alpha +j\bar{Y}_\alpha ) \in \n_\C $, where
$Z\in \a $ and $X_\alpha, Y_\alpha \in \g^\alpha $. Then,
\[ \la U,V \ra = \sum_{\alpha\in \Sigma^+} \Re j \kappa(Z,X_\alpha +j Y_\alpha) -
\Re j \kappa(Z,\bar{X}_\alpha +j\bar{Y}_\alpha) = 0, \]
since $\g^\alpha\perp_\kappa\g^\beta$ for $\alpha+\beta\neq 0$. 
\\
\par Having established (\ref{100}) the symplectic action of $T$ 
on $M_a$ will be Hamiltonian with moment map $\Phi$ if 
$\iota(\tilde{Z})\omega = d\Phi_Z$ holds for all $Z\in \t$. 
Fix $b\in M_a$ and $Y\in \k$. With the identification 
(\ref{ident}) we then we compute
\begin{eqnarray*}
d\Phi_Z(b)(\tilde{Y}_b) & = &  \frac{d}{dt} \Big|_{t=0} \Phi_Z (\tilde{b}(\exp(tY).b) \\
& = & \frac{d}{dt} \Big|_{t=0} \la i \Phi(\tilde{b}(\exp(tY)b), Z \ra \\
& = & \frac{d}{dt} \Big|_{t=0} \la i\Phi(\tilde{b}(b\exp(t{\rm Ad}(b^{-1})Y)) , Z\ra \\
& = & \frac{d}{dt} \Big|_{t=0} \la i \Im \log \tilde{a} (b\exp(t{\rm Ad}(b^{-1})Y)), Z\ra \\
& = & \langle i\Im \pr_{\a_\C} ({\rm Ad}(b^{-1})Y),Z \rangle \\
& = & \langle \pr_{\a_\C} ({\rm Ad}(b^{-1})Y),Z \rangle .
\end{eqnarray*}
For the last equality we have used the fact that $\a\perp \k$ with respect 
to $\la, \ra$ (cf. Lemma 3.2). 
\par On the other hand,
\begin{eqnarray*}
(\iota(\tilde{Z})\omega)_b (\tilde{Y}_b) & = & \omega_b(\tilde{Z}_b,\tilde{Y}_b) =  
\langle \pr_{\k_\C}({\rm Ad}(b^{-1})Z),{\rm Ad}(b^{-1})Y \rangle \\
& = & \langle {\rm Ad}(b^{-1})Z, \pr_{\b_\C}({\rm Ad}(b^{-1})Y) \rangle \\
&  = & \langle Z, \pr_{\b_\C}({\rm Ad}(b^{-1})Y) \rangle \\
& = & \langle Z, \pr_{\a_\C}({\rm Ad}(b^{-1})Y) \rangle .
\end{eqnarray*}
The last two equations hold because ${\rm Ad}(b^{-1})Z \in Z+\n_\C $, and 
$\la \t,\n_\C\ra=\{0\}$. 
\end{proof}

\subsection{Symplectic proof of the complex convexity theorem}

As $M_a$ is compact and the action of $T$ on $(M_a, \omega)$ is 
Hamiltonian, the Atiyah-Guillemin-Sternberg convexity theorem \cite{A, GS}  asserts 
$$\Phi(M_a)={\rm conv} (\Phi ({\rm Fix}(M_a)))\ . $$ 
In this formula ${\rm conv}(\cdot)$ denotes the convex hull of $(\cdot)$ 
and ${\rm Fix}(M_a)$ stands for the $T$-fixed points in $M_a$. 
Standard structure theory implies that ${\rm Fix}(M_a)={\mathcal W}.a$. 
We have thus proved:

\begin{theorem}\label{convthm} Let $G$ be a complex semisimple Lie group 
and $X\in \Omega$. Then 
\[ \Im(\log \tilde{a}(K\exp(iX))) = {\rm conv}(\mathcal{W}.X)\ .  \]
\end{theorem}

\section{The complex convexity theorem for non-complex groups}

For real groups the totally real $K$-orbits are no longer symplectic manifolds. However, they can be viewed as fixed point sets of an involution $\tau$ on the compact symplectic manifold $M_a$ as introduced in Chapter 3. We will define $\tau$ and show that it is compatible with the action of $T$ in a way that the symplectic convexity theorem from \cite{KO} can be applied.

Let $\g_0$ be a non-compact real form of the complex Lie algebra $\g$. It is no loss 
of generality if we assume that $\g_0$ is $\theta$-invariant. 
With $\k_0=\g_0\cap \k$ and $\p_0=\g_0\cap \p$ we then obtain a 
Cartan decomposition $\g_0=\k_0+\p_0$ of $\g_0$. We fix a maximal abelian subalgebra $\a_0$ of $\p_0$ 
which is contained in $\a$. Write $\Sigma_0=\Sigma(\g_0, \a_0)$ for the corresponding 
restricted root system and set 
\[ \Omega_0 = \{ X\in \a_0 : \abs{\alpha(X)}<
\frac{\pi}{2} \ \forall \alpha \in \Sigma_0 \}\ . \]
As $\Sigma_0=\Sigma|_{\a_0}\backslash  \{0\}$ we record that 
\begin{equation}
\Omega_0\subseteq \Omega \quad \hbox{and} \quad \Omega_0=\Omega\cap \a_0\ .
\end{equation}
It is no loss of generality to assume that 
$\Sigma_0^+=\Sigma^+|_{\a_0}\backslash\{0\}$ defines a positive system 
of $\Sigma_0^+$. We form the nilpotent 
Lie algebra $\n_0=\bigoplus_{\alpha\in \Sigma_0^+} \g_0^\alpha$ and record
that $\n_0=\g_0\cap\n$.

\par The analytic subgroups of $G$ with Lie algebras $\g_0,\k_0,\a_0$ and $\n_0$
will be denoted by $G_0, K_0, A_0$ and $N_0$.  If $\Xi_0=G_0\exp(i\Omega_0)(K_0)_\C / (K_0)_\C$
denotes the crown of $G_0/K_0$, then (4.1) yields a holomorphic $G_0$-equivariant 
embedding 
\begin{equation}
\Xi_0\to \Xi\ .
\end{equation}

As described at the end of Subsection 2.1, there exists a map
$\log \tilde{a}_0: \tilde{\Xi}_0=(N_0)_\C (A_0)_\C (K_0)_\C \to (\a_0)_\C$
with $\log \tilde{a}_0({\bf 1})=0$. Note that $\log \tilde{a}_0=\log \tilde{a}|_{\tilde{\Xi}_0}$.

\par Let $\sigma$ denote the Cartan involution on $\g_0$. We also write $\sigma$ for the 
doubly complex linear extension of $\sigma$ to $\g_\C$.  Likewise 
$\theta$ also stands for the complex linear extension of $\theta$ to 
$\g_\C$. We will be interested in the involution $\tau=\theta \circ \sigma=\sigma \circ \theta$ on $\g_\C$ (which is the complex linear extension of the complex conjugation on $\g$ with respect to $\g_0$). \\
All these involutions on $\g_\C$ can be lifted to involutions on $G_\C / K_\C $ and on $B_\C$, and we use the same letters to denote the lifts. \\

Recall the definition of the symplectic manifold $M_a=\tilde{b}(Ka)$ from Chapter 3. 
Notice that $M_a$ is $\tau$-invariant. The connected component of the $\tau$-fixed point set 
which contains $a$ is given by $Q_a=\tilde{b}(K_{0}a)$. 
We also have a Hamiltonian action by the torus $T_0 =\exp(j\a_0)$ on $M_a$. \\
The following lemma describes certain compatibility properties of the actions of $T_0$ and $\tau$ 
on $M_a$.

\begin{lemma}\label{compatibility}
Consider the Hamiltonian torus action of $T_0=\exp(j\a_0)$ on $M_a$ with momentum map $\Phi:M_a \rightarrow \t_0^*$. 
Then the following assertions hold: 
\begin{enumerate}
\item $t \circ \tau = \tau \circ t^{-1} \quad $ for all $t\in T_0$. \label{anticomm} \\
\item $\Phi \circ \tau = \Phi $.  \label{invariance} \\
\item $Q_a $ is a Lagrangian submanifold of $M_a$. \label{Lagr}
\end{enumerate}
\end{lemma}
\begin{proof}
(\ref{anticomm}) : \quad For $b=\tilde{b}(ka)\in M_a$ and $t\in T_0$,
\[ t.\tau(b)=t. \tilde{b}(\tau(k)a)=\tilde{b}(t\tau(k)a)=\tilde{b}(\tau(t^{-1}k)a)=\tau(t^{-1}.b). \]
(\ref{invariance}) : We need to introduce some additional notation. \\
Let
\[ \n^+ := \bigoplus_{\alpha\in \Sigma_0^+} \g^\alpha \subseteq \n , \quad
 \n^- := \bigoplus_{\alpha\in \Sigma_0^+} \g^{-\alpha} , \quad \mbox{and}
\quad \n^0 := \bigoplus_{\alpha\in \Sigma^+ \backslash \Sigma_0} \g^\alpha . \]
We denote by $N_\C^+$, $N_\C^-$ and $N_\C^0$ the analytic subgroups of $G_\C$ with Lie algebras $\n_\C^+$, $\n_\C^-$ and 
$\n_\C^0$, respectively. Notice that $\n=\n^0+\n^+$ and therefore $N_\C =N_\C^+ N_\C^0$.
It is important to observe that $\tau(N_\C^+)=N_\C^+$ but $\tau (N_\C^0)\cap N_\C^0=\{{\bf 1}\}$.  

\par Write  $\tilde{\t}$ for a $\tau$-invariant complement  of $\a_0 +i\a_0$ in $\a_\C$. 
Then $\a_\C=\a_0+i\a_0+\tilde{\t}$. \\
Let now $x=\tilde ka$ for some $\tilde k\in K$. Then 
$x, \tau(x)\in N_\C A_\C K_\C $ and  
\begin{equation}\label{Eins} x=n_+ n_0  b t k , \quad \tau(x)=n'_+ n'_0 b' t' k' 
\end{equation}
with elements $n_+ , n'_+ \in N_\C^+, \ n_0, n'_0 \in N_\C^0, \ b, b' \in \exp(\a_0+i\a_0), 
\ t, t' \in \exp(\tilde{\t})$ and $k, k' \in K_\C $. 
\par Clearly, (\ref{invariance}) will be proved if we can show that $b^2=(b')^2$ (which forces 
$b=b'$ by the comments at the end of Subsection 2.1). This will be established in the sequel. 

\par It follows from (\ref{Eins}) that 
\[ \tau(x) = n'_+ n'_0 b' t' k' = \tau(n_+) \tau(n_0)  b t^{-1} \tau(k) . \]
Since $\tau$ leaves $K_\C$ invariant and $\theta$ fixes each element of $K_\C$, we obtain
\begin{eqnarray}
\tau(x) \theta(\tau(x))^{-1} & = & n'_+ n'_0  (b')^2 \theta(n'_0)^{-1} \theta(n'_+)^{-1} \\
& = & \tau(n_+)\tau(n_0)b^2  \theta(\tau(n_0))^{-1} \theta(\tau(n_+))^{-1} . 
\end{eqnarray}
Notice that $n'_0 (b')^2  \theta(n'_0)^{-1}$ and $\tau(n_0)b^2 
\theta(\tau(n_0))^{-1}$ belong to the reductive group $Z_{G_\C}(A_0)$, and recall that $\tau(N_\C^+)=N_\C^+$
 and $\theta(N_\C^+)=N_\C^-$. Hence (4.4-5) combined with 
the Bruhat decomposition of $G_\C$ with respect to the parabolic subgroup $Z_{G_\C}(A_0) N_\C^+$ 
forces that $n'_+ = \tau(n_+)$. But then we have 
\[ n'_0(b')^2 \theta(n'_0)^{-1} = \tau(n_0)b^2  \theta(\tau(n_0))^{-1} \]
in $Z_{G_\C}(A_0)$. 
The components of $A_{0,\C}$, the center of $Z_{G_\C}(A_0)$, on both sides must coincide, therefore 
\[ (b')^2 = b^2 . \]
(\ref{Lagr}) : \quad Consider $U,V\in \k_0, k_0\in K_0$ and $b=\tilde{b}(k_0 a)\in Q_a$. From the formula (\ref{syform}) for the symplectic form $\omega$ on $M_a$ we get
\[ \omega_b(\tilde{U}_b,\tilde{V}_b)=\langle \pr_{\k_\C}({\rm Ad}(b^{-1})U),{\rm Ad}(b^{-1})V \rangle . \]
Now, both $\pr_{\k_\C}({\rm Ad}(b^{-1})U)$ and ${\rm Ad}(b^{-1})V$ lie in $\g_0 +i\g_0$. But for general elements $X_1, X_2, Y_1, Y_2 \in \g_0$ we have
\begin{eqnarray*}
\langle X_1+iX_2,Y_1+iY_2 \rangle & = & \Re \kappa(X_1+jX_2,Y_1+jY_2)-\Re \kappa (
\bar{X_1}+j\bar{X_2},\bar{Y_1}+j\bar{Y_2}) \\
& = & \Re \kappa(X_1,Y_1)-\Re \kappa(X_2,Y_2)-\Re \kappa(\bar{X_1},\bar{Y_1})+\Re \kappa(\bar{X_2},\bar{Y_2}) \\
&   & + \Re \kappa(X_1,jY_2)+\Re \kappa(jX_2,Y_1)-\Re \kappa(\bar{X_1},j\bar{Y_2})-\Re \kappa(j\bar{X_2},\bar{Y_1}) \\
& = & 0
\end{eqnarray*}
The last equality is due to the invariance of $\kappa$ and the fact that $X_1, X_2, Y_1, Y_2 \in \g_0$.
This shows that $\omega_b(\tilde{U}_b,\tilde{V}_b)=0$, i.e. $T_{b}(Q_a)$ is isotropic.
\end{proof}

We recall the following symplectic convexity theorem \cite{KO}.
\begin{theorem}
Let $M$ be a compact connected symplectic manifold with Hamiltonian torus action 
$T\times M\to M$ and momentum map $\Phi : M\to \t^*$. In addition, let $\tau : M\to M$ be an involutive diffeomorphism with fixed point set $Q$ such that 
\begin{enumerate}
\item $t\circ \tau=\tau\circ t^{-1}$ for all $t\in T$. \label{anticomm1} \\
\item $\Phi\circ \tau=\Phi$. \label{invariance1} \\
\item $Q$ is a Lagrangian submanifold of $M$. \label{Lagr1}
\end{enumerate} 
Denote the $T$-fixed subsets of $M$ and $Q$ by ${\rm Fix}(M)$ and ${\rm Fix}(Q)$, respectively.
Then, 
\[\Phi(Q)=\Phi(M)=\rm{conv}(\Phi({\rm Fix}(M)))=\rm{conv}(\Phi({\rm Fix}(Q))). \]
Moreover, the same assertions hold if $Q$ is replaced with
any of its connected components.
\end{theorem}

With this result at hand we are now able to prove the complex convexity result for non complex groups.
Write ${\mathcal W}_0$ for the Weyl group of $\Sigma_0$. 
\begin{theorem}
Let $G_0$ be a non-compact connected semisimple Lie group with Lie algebra $\g_0$. Fix an element $X\in \Omega_0$. Then
\begin{equation}\label{gleichung}
\Im \log \tilde{a}_0(K_0 \exp(iX))=\rm{conv}(\mathcal{W}_0.X)\ .
\end{equation}
\end{theorem}
\begin{proof}
Define $a=\exp(iX)$.  The left hand side in equality (\ref{gleichung}) coincides with $\Phi(Q_a)$ 
where $\Phi=\Im \circ \log \circ \tilde{a}$ is the momentum map on $M_a$. 
Lemma \ref{compatibility} says that conditions (\ref{anticomm1})-(\ref{Lagr1}) in Theorem 4.2 are satisfied. Therefore,
\[ \Im \log \tilde{a}_0(K_0 \exp(iX))=\Phi(Q_a)=\rm{conv}(\Phi({\rm Fix}(Q_a))). \]
Standard structure theory shows that ${\rm Fix}(Q_a)=\tilde{b}(\mathcal{W}_0\exp(iX))$.
This implies $\Phi({\rm Fix}(Q_a))=\mathcal{W}_0.X$, and finishes the proof.
\end{proof}

\bibliographystyle{amsplain}

\end{document}